\documentclass{amsart}
\usepackage{amsthm}
\usepackage{amsmath}
\usepackage{amsfonts}
\usepackage{amssymb}
\usepackage{graphicx}
\usepackage{float}
\usepackage{rotating}
\restylefloat{figure}
\usepackage{mathrsfs}
\usepackage{fancyhdr}
\usepackage[breaklinks]{hyperref}
\usepackage[numbers,sort&compress]{natbib}
\theoremstyle{plain}
\newtheorem{lem}{Lemma}
\newtheorem{prop}[lem]{Proposition}
\newtheorem{thm}[lem]{Theorem}
\newtheorem{cor}{Corollary}

\theoremstyle{definition}

\newtheorem{rem}[lem]{Remark}

\newcommand{\R}{\mathbb R}

\newcommand{\N}{\mathbb N}
\renewcommand{\d}{\,d}

\newcommand{\eps}{\varepsilon}
\renewcommand{\phi}{\varphi}
\newcommand{\norm}[1]{\Vert#1\Vert}

\renewcommand{\autoref}[1]{\text{Eq.}~\eqref{#1}}
\newcommand{\bea}{\begin{eqnarray}}
\newcommand{\eea}{\end{eqnarray}}
\newcommand{\beq}{\begin{equation}}
\newcommand{\eeq}{\end{equation}}
\begin{document}
\begin{center}
  {\small\rm Original publication in J. Math. Anal. Appl. 408(2) 513–-524 (2013)}
\end{center}
\title{A new model for electrostatic MEMS with two free boundaries}
\author{Martin Kohlmann}
\address{Dr.\ Martin Kohlmann, Goerdelerstraße 36, 38228 Salzgitter, Germany}
\email{martin\_kohlmann@web.de}
\date{\today}
\begin{abstract} A moving boundary problem with two free boundaries modeling a two-dimensional idealized MEMS device with pull-in instability is discussed. We use a fixed point argument to show that the model possesses stationary solutions for small source voltages. We also give a rigorous evidence that solutions of the model converge towards solutions of the associated small aspect ratio equation in the vanishing aspect ratio limit.
\end{abstract}
\maketitle
\tableofcontents
\section{Introduction and main results}\label{sec_intro}
Micro-Electro Mechanical Systems (MEMS) are a technology that relies on miniaturized mechanical elements and electro-mechanical devices obtained from microfabrication. Depending on their type of use, MEMS have various applications, e.g., as microsensors or microactuators. Over the past several decades researchers and developers have demonstrated miscellaneous types of MEMS for almost every possible sensing and actuating modality as temperature, pressure, inertial forces, chemical species, magnetic fields or radiation. In particular, accelerometers have become key components of commercial systems as smart phones, navigation systems, air bags and of health-enabling technologies. Due to the wide range of applications and their practical importance, the technology of MEMS has eventually created an interdisciplinary area of research coined by joint work of engineers, physicists and mathematicians.

A simple idealized mathematical model for a modern MEMS device consists of a rigid ground plate and a thin, deformable and elastic membrane that is suspended above the rigid plate and is held fixed along its boundary. A voltage difference between the ground plate and the free membrane induces an electric field in the space between the two components of the device and causes a Coulomb force on the membrane and hence a mechanical deformation. We refer the reader to \cite{P01,P0102,P03} for a detailed presentation of the modeling assumptions leading to the governing equations for this type of model.

\begin{figure}[H]
\begin{center}
\includegraphics[width=8cm]{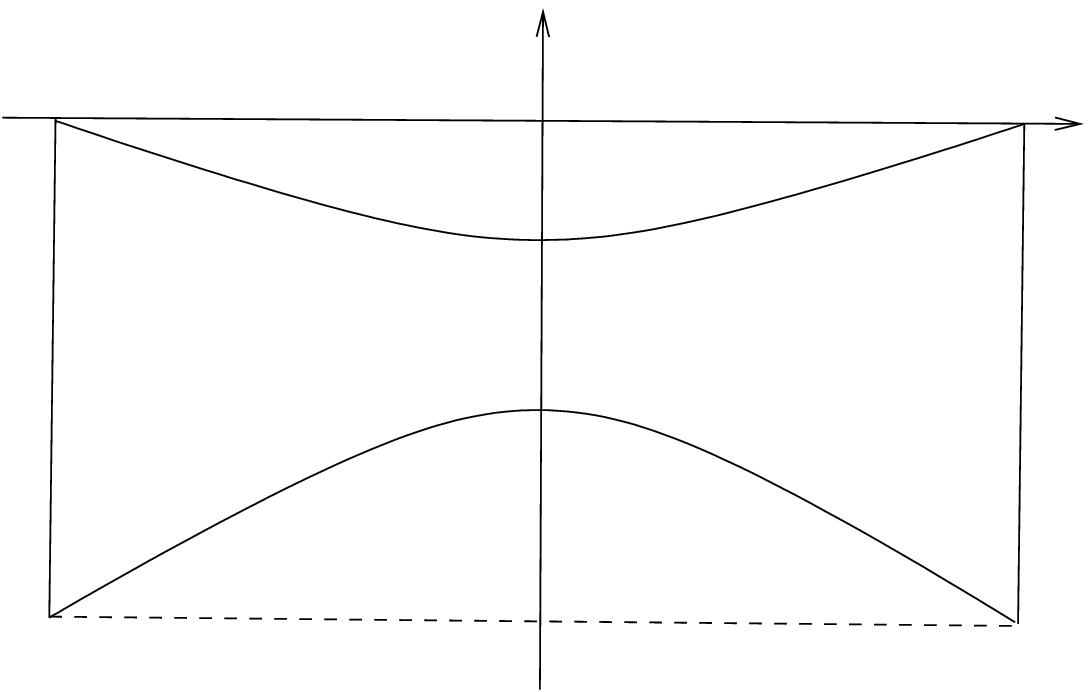}
\end{center}
\caption{An idealized model for an electrostatic MEMS device with two free boundaries.}
\label{figabstr}
\end{figure}

In the present paper, we intend to discuss an enhanced model for electrostatic MEMS consisting also of two components, but both of them are free elastic membranes. Pick $u,v\in W_\infty^2(-1,1)$, $0\geq u>v\geq-1$, let
$$\Omega_{u,v}=\{(x,z)\in(-1,1)\times\R;\,v(x)<z<u(x)\}$$
and denote the two horizontal boundary components of $\Omega_{u,v}$ by $\Gamma_u=\{z=u(x)\}$ and $\Gamma_v=\{z=v(x)\}$. The case when no source voltage is applied is modeled by $\Gamma_{0}$ and $\Gamma_{-1}$. For non-zero source voltage, the one-dimensional displacement of the membranes from $\Gamma_0$ and $\Gamma_{-1}$ respectively is modeled by the functions $u$ and $v$; see Fig.~\ref{figabstr}.
\begin{picture}(0,0)
\put(135,136){$\Omega_{u,v}$}
\put(110,111){$\Gamma_v$}
\put(110,159){$\Gamma_u$}
\put(190,183){\vector(0,-1){23}}
\put(190,78){\vector(0,1){40}}
\put(195,170){$u(x)$}
\put(195,88){$v(x)+1$}
\put(290,181){$z=0$}
\put(290,76){$z=-1$}
\put(260,186){$+1$}
\put(60,186){$-1$}
\put(172,203){$z$}
\put(277,173){$x$}
\end{picture}
Let $\partial_x=\frac{\partial}{\partial x}$, $\partial_z=\frac{\partial}{\partial z}$, $\nabla_\eps=(\eps\partial_x,\partial_z)$ and $\Delta_\eps=\eps^2\partial_x^2+\partial_z^2$. Our problem reads
\begin{align}
\Delta_\eps\phi & = 0, & \text{in }\Omega_{u,v},\label{originalproblem1}\\
\phi & = \frac{z-v}{u-v}, & \text{on }\partial\Omega_{u,v},\label{originalproblem2}\\
\partial_x^2u & = \lambda|\nabla_\eps\phi|^2,&\text{on }\Gamma_u,\label{originalproblem3}\\
\partial_x^2v & = - \mu|\nabla_\eps\phi|^2,&\text{on }\Gamma_v,\label{originalproblem4}\\
u(\pm 1) & = 0,\label{originalproblem5}\\
v(\pm 1) & =-1.\label{originalproblem6}
\end{align}
The function $\phi$ models the electrostatic potential in the region $\Omega_{u,v}$ and \eqref{originalproblem1} is the Laplace equation. The parameter $\eps>0$ is the \emph{aspect ratio} of the device comparing gap size to device length. The boundary condition \eqref{originalproblem2} stems from the assumption that the potential is equal to one on the upper membrane, vanishes on the lower membrane and is linear on the boundary components $\overline{\Omega_{u,v}}\cap\{x=\pm 1\}$. The parameters $\lambda,\mu>0$ in \autoref{originalproblem3} and \autoref{originalproblem4} are inversely proportional to the surface tension coefficients of the membranes. Precisely, the left-hand sides of \eqref{originalproblem3} and \eqref{originalproblem4} correspond to the curvature of $\Gamma_u$ and $\Gamma_v$ respectively, and the right-hand side is the energy density of the electric field inside $\Omega_{u,v}$. The minus sign in \eqref{originalproblem4} is a consequence of the fact that the membranes should attract each other. Note that $\lambda$ and $\mu$ depend quadratically on the source voltage of the device. Finally, we need the conditions \eqref{originalproblem5} and \eqref{originalproblem6} as both membranes should be held fixed on the boundary of $\Omega_{u,v}$. Observe that \eqref{originalproblem1}--\eqref{originalproblem6} is a free boundary problem since the domain $\Omega_{u,v}$ and its boundary components $\Gamma_u$ and $\Gamma_v$ are a priori unknown and depend on the solution $(u,v,\phi)$.

A major factor limiting the effectiveness of the idealized MEMS device is the so-called \emph{pull-in instability}: If the source voltage increases beyond a certain threshold value, the membranes $\Gamma_u$ and $\Gamma_v$ will come closer and closer and finally touch. In this case, there is no steady state of the MEMS but structure collapse and failure instead as the mechanical forces can no longer balance the electrostatic forces. To guarantee the functionality of a MEMS device it is of particular importance to know precisely the value of the pull-in voltage. In this context, \emph{small aspect ratio models} have been studied with regularity in recent years; we refer the reader to \cite{EGG10,P01,P03} for an overview and to \cite{BGP00,FMP03,GG0607,GG08,G08a,G08b,P0102,PD05} for further details. The recent references \cite{BL11,BP11,LW12,W12} relate generally to consequences of the small aspect ratio limit. Formally sending $\eps\to 0$ in \eqref{originalproblem1}--\eqref{originalproblem6} yields the solution
\beq\label{phi0}\phi_0(x,z)=\frac{z-v}{u-v}\eeq
and the coupled system
\begin{align}
u''(x)=\frac{\lambda}{(u-v)^2},&\quad u(\pm 1)=0,\label{usar}\\
v''(x)=-\frac{\mu}{(u-v)^2},&\quad v(\pm 1)=-1.\label{vsar}
\end{align}
\autoref{usar} and \autoref{vsar} will be referred to as the small aspect ratio equations of the model \eqref{originalproblem1}--\eqref{originalproblem6}. The pull-in instability corresponds to the singularity of the right-hand sides of \eqref{usar} and \eqref{vsar} obtained for $u(x)=v(x)$.

The small aspect ratio model of the MEMS model with only one free boundary on $(-1/2,1/2)$ reads
\beq\label{sareqtrsf}w''(x)=\frac{\Lambda}{(w+1)^2},\quad w(\pm 1/2)=0,\eeq
and has been widely studied in recent years. In \cite{BGP00,PT00} it is explained that \autoref{sareqtrsf} is analytically solvable and that its solutions $w(x;\Lambda)$ satisfy the implicit formula
$$\sqrt{\frac{(w+1)(w+1-\Lambda/E)}{2E}}+\frac{\Lambda}{E\sqrt{2E}}\tanh^{-1}\sqrt{\frac{w+1-\Lambda/E}{w+1}}=x,$$
where $E$ as a function of $\Lambda$ is obtained from
$$\sqrt{\frac{1-{\Lambda}/{E}}{2E}}+\frac{\Lambda}{E\sqrt{2E}}\tanh^{-1}\sqrt{1-\Lambda/E}=1/2.$$
It is shown numerically and by bifurcation arguments that there exists a threshold $\Lambda_*$ such that, for $0<\Lambda<\Lambda_*$, \autoref{sareqtrsf} has two solutions and for $\Lambda>\Lambda_*$, there exists no solution. One of the solutions in the case $0<\Lambda<\Lambda_*$ is stable under perturbations while the other one is instable. Since $\Lambda$ is quadratically proportional to the source voltage, the threshold $\Lambda_*$ corresponds to the pull-in voltage of the device. Moreover, an asymptotic expansion for $w$ is computed. The more general narrow-gap model
$$\alpha^2\partial_t^2w+\partial_tw-\partial_x^2w=-\frac{\Lambda}{(w+1)^2},\quad (t,x)\in(0,\infty)\times(-1,1),$$
has been subject of \cite{FMP03,G10,KLNT11} in the hyperbolic case $\alpha>0$ and of \cite{FMPS0607,GG08,G08a,G08b,H11,P0102} in the parabolic case $\alpha=0$.

For $v\equiv -1$ and $\mu=0$, we recover the MEMS model with a fixed ground plate from \eqref{originalproblem1}--\eqref{originalproblem6} and the problem \eqref{originalproblem1}--\eqref{originalproblem6} reduces to the model presented in \cite{WL12} where the authors prove the existence of a solution $(u,\phi)$ in suitable Sobolev and H\"older spaces for small values of the parameter $\lambda$, i.e., for small voltage values. In addition, a justification of the small aspect ratio model \eqref{sareqtrsf} is given by showing that steady state solutions of the free boundary problem converge towards stationary solutions of the small aspect ratio equation when the aspect ratio $\eps$ of the device tends to zero.

The plan of the paper at hand is first to obtain an existence result for the problem \eqref{originalproblem1}--\eqref{originalproblem6} with two free boundaries. We prove that for sufficiently small values of $\lambda,\mu$, there exists a solution $(u,v,\phi)$ to \eqref{originalproblem1}--\eqref{originalproblem6} in appropriate function spaces.
\thm\label{thm1} Fix some $\eps\in(0,1)$, $r_0\in(0,\frac{2}{3})$ and $\alpha\in[0,1)$. There is a number $a_0=a_0(r_0)>0$ so that the problem \eqref{originalproblem1}--\eqref{originalproblem6} admits for each $(\lambda,\mu)\in(0,a_0]^2$ a solution $(u,v,\phi)$ such that
\begin{align}
u,v&\in C^{2+\alpha}([-1,1]),\nonumber\\
\phi&\in W_2^2(\Omega_{u,v})\cap C(\overline{\Omega_{u,v}})\cap C^{2+\alpha}(\Omega_{u,v}\cup\Gamma_u\cup\Gamma_v),\nonumber
\end{align}
the functions $u,v$ and $x\mapsto\phi(x,z)$, for fixed $z$, are even and $u$ and $-v$ are convex and satisfy $|\partial_x^2u(x)|,|\partial_x^2v(x)|\leq r_0$.
Moreover, there is $\kappa_0\in(0,\tfrac{1}{3})$, independent of $\eps$, such that
$$0\geq u(x)\geq -\tfrac{1}{3}+\kappa_0\quad\text{and}\quad-\tfrac{1}{3}-2\kappa_0\geq v(x)\geq -1$$
and
$$\norm{u}_{W_\infty^2(-1,1)},\norm{v}_{W_\infty^2(-1,1)}\leq 3.$$
\endthm\rm
Our second result shows that for any family $\{(u_{\eps},v_{\eps},\phi_{\eps})\}_{{\eps}\in(0,1)}$ of solutions, obtained under the assumptions of Theorem~\ref{thm1}, there is a null sequence $\{\eps_k\}_{k\in\N}$ such that convergence $(u_{\eps_k},v_{\eps_k},\phi_{\eps_k})\to(u_0,v_0,\phi_0)$, $k\to\infty$, to a solution of the small aspect ratio model \eqref{phi0}--\eqref{vsar} associated with \eqref{originalproblem1}--\eqref{originalproblem6} holds. Let $\mathbf 1_{A}$ denote the characteristic function of the set $A\subset\R^2$.
\thm\label{thm2} Let $\{(u_\eps,v_\eps,\phi_\eps)\}_{\eps\in(0,1)}$ be a family of solutions to \eqref{originalproblem1}--\eqref{originalproblem6} as specified in Theorem~\ref{thm1}. Then there exists a null sequence $\{\eps_k\}_{k\in\N}$ and solutions $u_0,v_0$ to the small aspect ratio equations \eqref{usar} and \eqref{vsar} such that
$$u_{\eps_k}\to u_0\quad\text{and}\quad v_{\eps_k}\to v_0\quad\text{in }W_\infty^1(-1,1)$$
and
$$\phi_{\eps_k}\mathbf 1_{\Omega_{u_{\eps_k},v_{\eps_k}}}\to\phi_0\mathbf 1_{\Omega_{u_0,v_0}}\quad\text{in }L_2\left((-1,1)\times(0,1)\right),$$
as $k\to\infty$, where $\phi_0=\frac{z-u_0}{u_0-v_0}$ is the potential in \autoref{phi0}.
\endthm\rm
\section{Solvability beyond the pull-in instability}
\begin{figure}[H]
\begin{center}
\includegraphics[width=11cm]{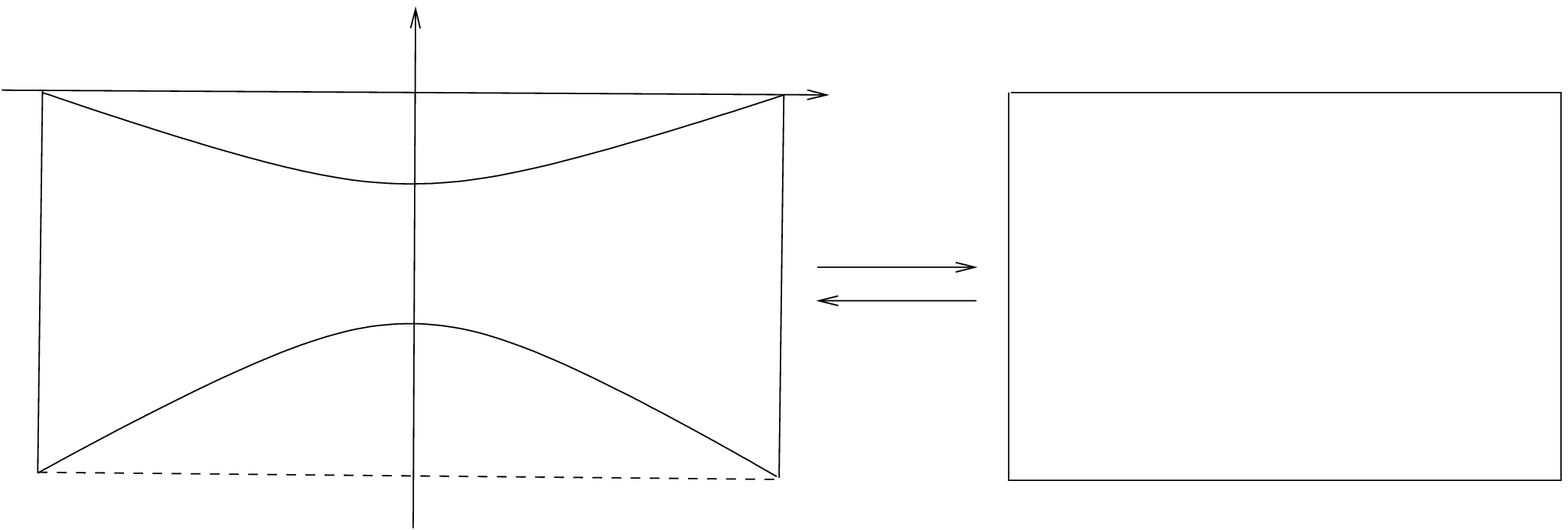}
\end{center}
\caption{Transformation onto a fixed reference domain.}\label{figtrsf}
\end{figure}

We first transform the problem \eqref{originalproblem1}--\eqref{originalproblem6} on the a priori unknown domain $\Omega_{u,v}$ to a corresponding problem on the fixed reference domain $\Omega:=(-1,1)\times(0,1)$; see Fig.~\ref{figtrsf}. Therefore, we introduce a transformation of coordinates $T=T_{u,v}\colon\overline{\Omega_{u,v}}\to\overline{\Omega}$ given by
\beq\label{defT}T(x,z)=(x',z')=\left(x,\frac{z-v(x)}{u(x)-v(x)}\right).\eeq
One easily checks that, for any $u,v\in W_\infty^2(-1,1)$, $T_{u,v}$ is a diffeomorphism $\overline{\Omega_{u,v}}\to\overline{\Omega}$ with the inverse
$$T^{-1}(x',z')=(x',z'(u(x')-v(x'))+v(x')).$$
\begin{picture}(0,0)
\put(70,105){$\Omega_{u,v}$}
\put(227,98){$\Omega=(-1,1)\times(0,1)$}
\put(182,112){$T_{u,v}$}
\put(182,85){$T_{u,v}^{-1}$}
\put(98,153){$z$}
\put(170,143){$x$}
\put(200,53){$x'=-1$}
\put(310,53){$x'=1$}
\put(255,65){$z'=0$}
\put(255,130){$z'=1$}
\end{picture}
Let $\widetilde{\nabla}_{u,v;\eps}=T_{u,v}\circ\nabla_\eps\circ T_{u,v}^{-1}$ and $\widetilde{\Delta}_{u,v;\eps}=T_{u,v}\circ\Delta_\eps\circ T_{u,v}^{-1}$ denote the transformed first and second order derivative operators on $\Omega$. We omit the tedious calculations leading to the explicit formulas
\begin{align}\widetilde{\Delta}_{u,v;\eps}w & = \eps^2w_{x'x'} - 2\eps^2w_{x'z'}\frac{z'(\tilde u'-\tilde v')+\tilde v'}{\tilde u-\tilde v}
+w_{z'z'}\frac{1+\eps^2[z'(\tilde u'-\tilde v')+\tilde v']^2}{(\tilde u-\tilde v)^2} \nonumber\\
& \quad +\eps^2w_{z'}\left(2\frac{\tilde u'-\tilde v'}{(\tilde u-\tilde v)^2}[z'(\tilde u'-\tilde v')+\tilde v']-\frac{z'(\tilde u''-\tilde v'')+\tilde v''}{\tilde u-\tilde v}\right)\label{deltatilde}
\end{align}
and
$$
\widetilde{\nabla}_{u,v;\eps}w=\begin{pmatrix}\eps w_{x'}\\0\end{pmatrix}+\frac{w_{z'}}{\tilde u-\tilde v}
\begin{pmatrix}\eps[z'(\tilde v'-\tilde u')-\tilde v']\\1\end{pmatrix}
$$
with $(\tilde u,\tilde v,\tilde\phi)=(u,v,\phi)\circ T_{u,v}^{-1}$. The transformed problem \eqref{originalproblem1}--\eqref{originalproblem6} on the fixed domain $\Omega$ then reads
\begin{align}
\widetilde{\Delta}_{u,v;\eps}\tilde\phi & = 0, & \text{in }\Omega,\label{Laplaceeq}\\
\tilde\phi & = z', & \text{on }\partial\Omega,\label{Laplacebc}\\
\partial_{x'}^2\tilde u & = \lambda|(\widetilde{\nabla}_{u,v;\eps}\tilde\phi)(\cdot,1)|^2,&\text{on }(-1,1)\label{sodu},\\
\partial_{x'}^2\tilde v & = -\mu|(\widetilde{\nabla}_{u,v;\eps}\tilde\phi)(\cdot,0)|^2,&\text{on }(-1,1)\label{sodv},\\
\tilde u(\pm 1) & = 0,\label{u(pm1)}\\
\tilde v(\pm 1) & =-1.\label{v(pm1)}
\end{align}
For $r_0\in(0,\frac{2}{3})$, let
\begin{align}
\mathcal C_1&:= \{u\in W_\infty^2(-1,1);\,u(x)=u(-x),\,u(\pm 1)=0,\,0\leq\partial_x^2u\leq r_0\}, \nonumber\\
\mathcal C_2&:= \{v\in W_\infty^2(-1,1);\,v(x)=v(-x),\,v(\pm 1)=-1,\,-r_0\leq\partial_x^2v\leq 0\}. \nonumber
\end{align}
The following lemma and Lemma 4 of \cite{WL12} provide some important properties of the sets $\mathcal C_1$ and $\mathcal C_2$.
\lem\label{lem_C1C2} The sets $\mathcal C_1,\mathcal C_2$ are convex, closed and bounded subsets of $W_q^2(-1,1)$ for any $1\leq q\leq\infty$, and there is $\kappa_0\in(0,\tfrac{1}{3})$ such that $0\geq u(x)\geq-\frac{1}{3}+\kappa_0$ and $-\tfrac{1}{3}-2\kappa_0\geq v(x)\geq -1$, for all $(u,v)\in\mathcal C_1\times\mathcal C_2$, $-1\leq x\leq 1$.
\endlem
\proof It follows immediately from the definitions that $\mathcal C_1,\mathcal C_2$ are convex and closed in $W_\infty^2(-1,1)$. Thus $\mathcal C_1,\mathcal C_2$ are weakly closed in $W_\infty^2(-1,1)$ and hence convex and closed in $W_q^2(-1,1)$ for any $1\leq q\leq\infty$. As shown in \cite{WL12}, any $u\in\mathcal C_1$ satisfies
\begin{align}
&0\geq u\geq -\frac{r_0}{2},\nonumber\\
&|\partial_xu|\leq 2r_0\label{boundu'}
\end{align}
on $(-1,1)$, and $\mathcal C_1$ is bounded in $W^2_q(-1,1)$ for any $1\leq q\leq\infty$. Pick $v\in\mathcal C_2$ and $x,y\in(-1,1)$. 
We integrate
%
%
$$\partial_xv(x) = \partial_xv(y) + \int_y^x\partial_x^2v(z)\d z$$
with respect to $y$ on $(-1,1)$ and use that $v(-1)=v(+1)$ to obtain that
\beq\label{boundv'}|\partial_xv|\leq 2r_0.\eeq
Since $v$ is concave and $v(\pm 1)=-1$, we also have $v\geq -1$. Assume that $v$ attains a maximum at $x_m\in[0,1)$. Then $\partial_xv(x_m)=0$ and
$$v(x)-v(x_m) = \int_{x_m}^x\partial_xv(y)\d y = \int_{x_m}^{x}(x-y)\partial_x^2v(y)\d y.$$
Using once again the boundary condition $v(1)=-1$, we observe that
$$-1-v(x_m) = \int_{x_m}^1(1-y)\partial_x^2v(y)\d y\geq -r_0.$$
Then
$$v\leq v(x_m)\leq r_0-1<-\frac{r_0}{2},$$
as $r_0\in(0,\tfrac{2}{3})$, and we infer that
\beq\label{bounduv}-1\leq v<u\leq 0\quad\text{on }[-1,1].\eeq
In particular, $\mathcal C_2$ is bounded in $W^2_\infty(-1,1)$ and hence in $W^2_q(-1,1)$ for any $1\leq q\leq\infty$. We let $\kappa_0=\tfrac{1}{2}(\tfrac{2}{3}-r_0)$ to complete the proof.
\endproof
Lemma~\ref{lem_C1C2} shows that, for $(u,v)\in\mathcal C_1\times\mathcal C_2$, $\Omega_{u,v}$ is as defined in the introduction. The next step will be to solve the subproblem \eqref{Laplaceeq}, \eqref{Laplacebc} for $(u,v)\in\mathcal C_1\times\mathcal C_2$ given.
\prop\label{prop_solphi} Given $(u,v)\in\mathcal C_1\times\mathcal C_2$, there is a unique solution $\tilde\phi=\tilde\phi_{u,v;\eps}$ to the problem \eqref{Laplaceeq}, \eqref{Laplacebc}. Moreover, $\tilde\phi$ is even with respect to $x'$,
\beq -1\leq z'(\tilde u(x')-\tilde v(x'))+\tilde v(x')\leq \tilde\phi(x',z')\leq 1,\quad (x',z')\in \Omega,\label{estimate_tildephi}\eeq
and
\beq\label{estimate_norm}\norm{\tilde\phi}_{W_2^2(\Omega)}\leq c_1\eeq
for some positive constant $c_1=c_1(r_0,\eps)$. Finally, for any $\eps>0$, the mapping
$$(u,v)\mapsto\tilde\phi_{u,v;\eps}\colon\mathcal C_1\times\mathcal C_2\to W_2^2(\Omega)$$
is continuous when $\mathcal C_1\times\mathcal C_2$ is equipped with the product topology of $W_2^2(-1,1)\times W_2^2(-1,1)$.
\endprop\rm
\proof We show that $-\widetilde\Delta_{u,v;\eps}$ is elliptic for $(u,v)\in\mathcal C_1\times\mathcal C_2$ by showing that its principal part
$$
A:=\begin{pmatrix}
\eps^2 & -\eps^2\frac{z'(\tilde u'-\tilde v')+\tilde v'}{\tilde u-\tilde v} \\
-\eps^2\frac{z'(\tilde u'-\tilde v')+\tilde v'}{\tilde u-\tilde v} & \frac{1+\eps^2[z'(\tilde u'-\tilde v')+\tilde v']^2}{(\tilde u-\tilde v)^2}
\end{pmatrix}
$$
is positive definite. Let
$$t=\eps^2+\frac{1+\eps^2[z'(\tilde u'-\tilde v')+\tilde v']^2}{(\tilde u-\tilde v)^2}\quad\text{and}\quad d=\frac{\eps^2}{(\tilde u-\tilde v)^2}$$
denote the trace and the determinant of $A$. Using the estimates \eqref{boundu'}--\eqref{bounduv}
it is easy to derive that
$$1+\eps^2\leq t\leq C_1\eps^2 + C_2\quad\text{and}\quad d\geq\eps^2,$$
where $C_1,C_2>0$ only depend on $r_0$. Since the eigenvalues $e_\pm$ of $A$ are given by
$$e_\pm = \frac{1}{2}\left(t\pm\sqrt{t^2-4d}\right)$$
we conclude that
$$e_+\geq e_-\geq\frac{e_+e_-}{e_++e_-}=\frac{d}{t}\geq\frac{\eps^2}{C_1\eps^2+C_2}>0.$$
Hence $-\widetilde\Delta_{u,v;\eps}$ is strictly elliptic and the ellipticity constant depends only on $r_0$ and $\eps$ and not on $(u,v)\in\mathcal C_1\times\mathcal C_2$.

Setting $\psi(x',z')=\tilde\phi(x',z')-z'$, the problem \eqref{Laplaceeq}, \eqref{Laplacebc} is equivalent to
\begin{align}
-\widetilde{\Delta}_{u,v;\eps}\psi & = f_{u,v;\eps} & \text{in }\Omega,\label{Laplaceeq2}\\
\psi & = 0 & \text{on }\partial\Omega,\label{Laplacebc2}
\end{align}
where
\beq f_{u,v;\eps} = \eps^2\left(2\frac{\tilde u'-\tilde v'}{(\tilde u-\tilde v)^2}[z'(\tilde u'-\tilde v')+\tilde v']-\frac{z'(\tilde u''-\tilde v'')+\tilde v''}{\tilde u-\tilde v}\right).\label{deff}\eeq
Clearly, $f\in L_\infty(\Omega)$ with $\norm{f_{u,v;\eps}}_{L_\infty(\Omega)} \leq C_3\eps^2$, $C_3=C_3(r_0)>0$.
Rewriting $\widetilde{\Delta}_{u,v;\eps}$ in divergence form,
\begin{align}
\widetilde{\Delta}_{u,v;\eps}w & = \partial_{x'} \left(\eps^2w_{x'}-\eps^2\frac{z'(\tilde u'-\tilde v')+\tilde v'}{\tilde u-\tilde v}w_{z'}\right) \nonumber\\
& \quad + \partial_{z'} \left(-\eps^2\frac{z'(\tilde u'-\tilde v')+\tilde v'}{\tilde u-\tilde v}w_{x'} + \frac{1+\eps^2[z'(\tilde u'-\tilde v')+\tilde v']^2}{(\tilde u-\tilde v)^2}w_{z'}\right) \nonumber\\
& \quad -\eps^2(\tilde u'-\tilde v')\frac{z'(\tilde u'-\tilde v')+\tilde v'}{(\tilde u-\tilde v)^2}w_{z'}+\eps^2\frac{\tilde u'-\tilde v'}{\tilde u-\tilde v}w_{x'},\label{Lapdivform}
\end{align}
we see that the coefficients of $\widetilde{\Delta}_{u,v;\eps}$ are uniformly bounded with respect to $(u,v)\in\mathcal C_1\times\mathcal C_2$. Applying Theorem 9.1 and Theorem 10.1 in \cite{LU68}, we conclude that there exists a unique solution $\psi=\psi_{u,v;\eps}\in\mathring{W}_2^2(\Omega)$ to \eqref{Laplaceeq2}, \eqref{Laplacebc2} such that
\beq\label{estimate_LU68}\norm{\psi}_{W_2^2(\Omega)}\leq C_4\left(\norm{\psi}_{L_2(\Omega)}+1\right)\eeq
with a constant $C_4$ depending on $r_0$ and $\eps$, but not on $(u,v)\in\mathcal C_1\times\mathcal C_2$. The function $\tilde\phi(x',z')=\psi(x',z')+z'$ clearly solves \eqref{Laplaceeq} and \eqref{Laplacebc}. Let $\mathbf 1$ be the constant function with value $1$ on $\overline\Omega$. Then $\widetilde{\Delta}_{u,v;\eps}\mathbf 1=\widetilde{\Delta}_{u,v;\eps}\tilde\phi=0$ in $\Omega$ and $\tilde\phi\leq\mathbf 1$ on $\partial\Omega$. By the maximum principle $\tilde\phi\leq\mathbf 1$ in $\overline\Omega$ and hence $\tilde\phi(x',z')\leq 1$ for all $(x',z')\in\overline\Omega$. Let $w:=z'(\tilde u-\tilde v)+\tilde v$. By \eqref{deltatilde}, $\widetilde{\Delta}_{u,v;\eps}w=0$ in $\Omega$ and $w\leq\tilde\phi$ on $\partial\Omega$. A further application of the maximum principle achieves $w\leq\tilde\phi$ in $\overline\Omega$ and this completes the proof of \eqref{estimate_tildephi}. Combining \eqref{estimate_tildephi} and \eqref{estimate_LU68} we immediately obtain $\norm{\tilde\phi}_{W_2^2(\Omega)}\leq C_5$, for some positive constant $C_5$ only depending on $r_0$ and $\eps$, and hence \eqref{estimate_norm}. To see that $\tilde\phi$ is even, we observe that, since $\tilde u$ and $\tilde v$ are even, the operator $\widetilde{\Delta}_{u,v;\eps}$ is invariant under the transformation $(x',z')\mapsto(-x',z')$. Hence $\tilde\phi(-x',z')$ is another solution of \eqref{Laplaceeq}, \eqref{Laplacebc} and, by uniqueness, $\tilde\phi(x',z')=\tilde\phi(-x',z')$, as was to be shown. That $(u,v)\mapsto\tilde\phi_{u,v;\eps}$ is continuous follows from analogous arguments as in the proof of Lemma 6 in \cite{WL12}.
\endproof
In the following lemma we obtain control of the gradient of $\tilde\phi$ on the horizontal boundary components.
\lem\label{lem_estimatesolphi} There is a constant $c_2>0$ depending only on $r_0\in(0,\frac{2}{3})$ and $\eps\in(0,1)$ such that the solution $\tilde\phi_{u,v;\eps}$ to \eqref{Laplaceeq}, \eqref{Laplacebc}, for $(u,v)\in\mathcal C_1\times C_2$, satisfies
$$\norm{(\partial_{z'}\tilde\phi_{u,v;\eps})(\cdot,0)}_{W_2^{1/2}(-1,1)}, \norm{(\partial_{z'}\tilde\phi_{u,v;\eps})(\cdot,1)}_{W_2^{1/2}(-1,1)}\leq c_2$$
and there is a constant $c_3>0$ such that
$$0\leq(\partial_{z'}\tilde\phi_{u,v;\eps})(x',0),(\partial_{z'}\tilde\phi_{u,v;\eps})(x',1)\leq 1+c_3\eps^2,\quad x'\in(-1,1).$$
\endlem
\proof Applying Theorem II-5.5 in \cite{N67} and \eqref{estimate_norm},
there is a positive constant $k_1$ only depending on $\Omega$ such that
$$\norm{(\partial_{z'}\tilde\phi_{u,v;\eps})(\cdot,0)}_{W_2^{1/2}(-1,1)},\norm{(\partial_{z'}\tilde\phi_{u,v;\eps})(\cdot,1)}_{W_2^{1/2}(-1,1)}\leq k_1\norm{\tilde\phi_{u,v;\eps}}_{W_2^2(\Omega)}\leq c_2.$$
Let $w_{\alpha}(z'):=(z')^{1+\alpha}$, $\alpha>0$. Then $\tilde\phi_{u,v;\eps}\geq w_\alpha$ on $\partial\Omega$, and using \eqref{boundu'}, \eqref{boundv'}, \eqref{bounduv}, $|\tilde u''|,|\tilde v''|\leq r_0$ and $z'\in(0,1)$, we infer that there is a constant $k_2>0$ such that for interior points we have
\begin{align}
\widetilde{\Delta}_{u,v;\eps}w_\alpha & = \alpha(1+\alpha)(z')^{\alpha-1}\frac{1+\eps^2[z'(\tilde u-\tilde v)+\tilde v]^2}{(\tilde u-\tilde v)^2}\nonumber\\
& \quad + \eps^2(1+\alpha)(z')^\alpha\left(2\frac{\tilde u'-\tilde v'}{(\tilde u-\tilde v)^2}[z'(\tilde u'-\tilde v')+\tilde v']-\frac{z'(\tilde u''-\tilde v'')+\tilde v''}{\tilde u-\tilde v}\right)\nonumber\\
&\geq \frac{(1+\alpha)(z')^{\alpha-1}}{(\tilde u-\tilde v)^2}\left(\alpha-\tfrac{3}{2}k_2r_0\eps^2\right).\nonumber
\end{align}
As $r_0<\frac{2}{3}$ we conclude that $-\widetilde\Delta_{u,v;\eps}w_{k_2\eps^2}\leq 0$ in $\Omega$ so that, by the weak maximum principle,
$$w_{k_2\eps^2} - \tilde\phi_{u,v;\eps} \leq \max_{\overline\Omega}\{w_{k_2\eps^2} - \tilde\phi_{u,v;\eps}\} = \max_{\partial\Omega}\{w_{k_2\eps^2} - \tilde\phi_{u,v;\eps}\} \leq 0,$$
i.e., $\tilde\phi_{u,v;\eps}\geq w_{k_2\eps^2}$ in $\overline \Omega$. For $z'\neq 1$ we observe
\begin{align}
\frac{1}{z'-1}(\tilde\phi_{u,v;\eps}(x',z')-\tilde\phi_{u,v;\eps}(x',1)) & = \frac{1}{z'-1}(\tilde\phi_{u,v;\eps}(x',z')-1) \nonumber\\
& \leq\frac{1}{z'-1}(w_{k_2\eps^2}(z')-w_{k_2\eps^2}(1)).\nonumber
\end{align}
Sending $z'\to 1$ yields $(\partial_{z'}\tilde\phi_{u,v;\eps})(x',1)\leq(\partial_{z'}w_{k_2\eps^2})(1)=1+k_2\eps^2$. Finally, as $\tilde\phi_{u,v;\eps}\leq 1$ in $\Omega$ and $\tilde\phi_{u,v;\eps}(x',1)=1$, we must have $(\partial_{z'}\tilde\phi_{u,v;\eps})(x',1)\geq 0$. Similarly, we obtain
\begin{align}
\frac{1}{z'}(\tilde\phi_{u,v;\eps}(x',z')-\tilde\phi_{u,v;\eps}(x',0)) & = \frac{1}{z'}\tilde\phi_{u,v;\eps}(x',z') \nonumber\\
& \geq\frac{1}{z'}(w_{k_2\eps^2}(z')-w_{k_2\eps^2}(0)),\nonumber
\end{align}
and hence $(\partial_{z'}\tilde\phi_{u,v;\eps})(x',0)\geq (\partial_{z'}w_{k_2\eps^2})(0)=0$. Let $\tilde w_\alpha(z')=w_{\alpha}(1-z')$. Then
\begin{align}
\widetilde{\Delta}_{u,v;\eps}\tilde w_\alpha \geq \frac{(1+\alpha)(1-z')^{\alpha-1}}{(\tilde u-\tilde v)^2}\left(\alpha-\tfrac{3}{2}k_2r_0\eps^2\right),\nonumber
\end{align}
and $-\widetilde\Delta\tilde w_{k_2\eps^2}\leq 0$. Again, by the weak maximum principle,
$$\tilde w_{k_2\eps^2}-(1-\tilde\phi_{u,v;\eps})\leq\max_{\overline\Omega}\{\tilde w_{k_2\eps^2}-(1-\tilde\phi_{u,v;\eps})\}
=\max_{\partial\Omega}\{\tilde w_{k_2\eps^2}-(1-\tilde\phi_{u,v;\eps})\}\leq 0,$$
i.e., $\tilde\phi_{u,v;\eps}\leq 1-\tilde w_{k_2\eps^2}$ in $\overline\Omega$. Hence
\begin{align}
\frac{1}{z'}(\tilde\phi_{u,v;\eps}(x',z')-\tilde\phi_{u,v;\eps}(x',0)) & = \frac{1}{z'}\tilde\phi_{u,v;\eps}(x',z') \nonumber\\
& \leq\frac{1}{z'}\left(1-(1-z')^{1+k_2\eps^2}\right)\nonumber
\end{align}
and L'Hospital's rule implies that $(\partial_{z'}\tilde\phi_{u,v;\eps})(x',0)\leq 1+k_2\eps^2$. This completes the proof of the lemma.
\endproof
Equations \eqref{sodu} and \eqref{sodv} with the boundary conditions \eqref{u(pm1)} and \eqref{v(pm1)} read
\begin{align}
\partial_{x'}^2\tilde u & = \lambda\frac{1+\eps^2(\partial_{x'}\tilde u)^2}{(\tilde u-\tilde v)^2}\left|(\partial_{z'}\tilde\phi)(\cdot,1)\right|^2, \label{sodu2}\\
\tilde u(\pm 1) & = 0 \nonumber
\end{align}
and
\begin{align}
\partial_{x'}^2\tilde v & = -\mu\frac{1+\eps^2(\partial_{x'}\tilde v)^2}{(\tilde u-\tilde v)^2}\left|(\partial_{z'}\tilde\phi)(\cdot,0)\right|^2, \label{sodv2}\\
\tilde v(\pm 1) & = -1.\nonumber
\end{align}
Let
\begin{align}
g_{u,v}&:=\frac{1+\eps^2(\partial_{x'}\tilde u)^2}{(\tilde u-\tilde v)^2}|(\partial_{z'}\tilde\phi)(\cdot,1)|^2,\nonumber\\
h_{u,v}&:=\frac{1+\eps^2(\partial_{x'}\tilde v)^2}{(\tilde u-\tilde v)^2}|(\partial_{z'}\tilde\phi)(\cdot,0)|^2.\nonumber
\end{align}
Since $g_{u,v},h_{u,v}\in L_\infty(-1,1)$ the second order problem
\begin{align}
w''(x')&=(\lambda g_{u,v}(x'),-\mu h_{u,v}(x')),\nonumber\\
w(\pm 1)&=(0,-1)\nonumber
\end{align}
possesses for any $\lambda,\mu>0$ a unique solution $w=S(u,v)\in W_\infty^2(-1,1)^2$. Our goal is to show that the map $S=(S_1,S_2)$ has a fixed point in $\mathcal C_1\times\mathcal C_2$ by means of Schauder's theorem.
\prop\label{prop2.6} Pick $r_0\in(0,\tfrac{2}{3})$ and $\eps\in(0,1)$ and let
$$a_0=a_0(r_0):=\frac{r_0(\frac{3}{2}r_0-1)^2}{(1+c_3)^2(1+4r_0^2)}.$$
The problem \eqref{Laplaceeq}--\eqref{v(pm1)} admits for each $(\lambda,\mu)\in(0,a_0]^2$ a solution
$$(\tilde u,\tilde v,\tilde\phi_{u,v;\eps})\in \mathcal C_1\times\mathcal C_2\times W_2^2(\Omega).$$
The function $\tilde\phi_{u,v;\eps}$ satisfies the properties stated in Proposition~\ref{prop_solphi} and $\tilde u,\tilde v$ satisfy the bounds established in Lemma~\ref{lem_C1C2}. Moreover, $\tilde u,\tilde v$ belong to $W_2^{2+\sigma}(-1,1)$ for any $\sigma\in[0,1/2)$.
\endprop\rm
\proof Recall that $\mathcal C_1,\mathcal C_2$ are convex, closed and bounded subsets of $W_2^2(-1,1)$.
Using the arguments in the proof of Lemma 8 of \cite{WL12}, we see that $S_i\colon\mathcal C_1\times\mathcal C_2\to W_2^{2+\sigma}(-1,1)$, $i=1,2$, is continuous for each $\sigma\in[0,1/2)$ and there is a positive constant $c_4$ depending on $r_0$, $\eps$ and $\sigma$ such that
$$\norm{S(u,v)}_{W_2^{2+\sigma}(-1,1)^2}\leq\max\{\lambda,\mu\} c_4,\quad\forall (u,v)\in\mathcal C_1\times\mathcal C_2.$$
Clearly, $S(u,v)$ is even, since $\tilde u,\tilde v$ and $x'\mapsto\tilde\phi$ are even, and $S_1(u,v)$ and $-S_2(u,v)$ are convex in view of \eqref{sodu2} and \eqref{sodv2}. Observe that, for all $(u,v)\in\mathcal C_1\times\mathcal C_2$,
$$0\leq\partial_{x'}^2S_1(u,v)=\lambda\frac{1+\eps^2(\partial_{x'}\tilde u)^2}{(\tilde u-\tilde v)^2}|(\partial_{z'}\tilde\phi)(\cdot,1)|^2\leq \frac{\lambda}{a_0}r_0.$$
Then $0\leq\partial_{x'}^2S_1(u,v)\leq r_0$ for $\lambda\in(0,a_0]$. Similarly, one has $0\geq\partial_{x'}^2S_2(u,v)\geq -r_0$ for $\mu\in(0,a_0]$. This shows that $S\colon\mathcal C_1\times\mathcal C_2\to\mathcal C_1\times\mathcal C_2$ for $(\lambda,\mu)\in(0,a_0]^2$. The embedding $W_{2}^{2+\sigma}(-1,1)\hookrightarrow W_2^2(-1,1)$ is compact and hence $S\colon\mathcal C_1\times\mathcal C_2\to\mathcal C_1\times\mathcal C_2$ is continuous and compact and has a fixed point $(\tilde u,\tilde v)\in\mathcal C_1\times\mathcal C_2$ in view of Schauder's theorem. It follows from our definitions that $\tilde u$ and $\tilde v$ enjoy the properties states in the proposition.
\endproof
We finally improve the regularity of $(\tilde u,\tilde v,\tilde\phi)$ and pull the solution back to $\Omega_{u,v}$. The following corollary completes the proof of Theorem~\ref{thm1}.
\cor Let $(\tilde u,\tilde v,\tilde\phi_{u,v;\eps})$ denote the solution to \eqref{Laplaceeq}--\eqref{v(pm1)} obtained in Proposition~\ref{prop2.6} and let $\alpha\in[0,1)$. Then $(u,v,\phi)=(\tilde u,\tilde v,\tilde\phi_{u,v;\eps})\circ T_{u,v}$ is a solution to \eqref{originalproblem1}--\eqref{originalproblem6} with regularity
\begin{align}
u,v&\in C^{2+\alpha}([-1,1]),\nonumber\\
\phi&\in W_2^2(\Omega_{u,v})\cap C(\overline{\Omega_{u,v}})\cap C^{2+\alpha}(\Omega_{u,v}\cup\Gamma_u\cup\Gamma_v).\nonumber
\end{align}
\endcor
\proof As $\Omega_{u,v}$ is a Lipschitz domain (cf.\ \cite[Ch.\ 6.2]{GT}) and satisfies the exterior cone condition (cf.\ \cite[p.\ 193]{GT}) on any point of $\partial\Omega_{u,v}$, and as $T_{u,v}\in C(\overline{\Omega_{u,v}},\overline\Omega)$, we can apply the arguments in the proof of Corollary 10 of \cite{WL12} to infer that
$$\phi\in W_2^2(\Omega_{u,v})\cap C(\overline{\Omega_{u,v}})$$
is an even solution with respect to $x$ to \eqref{originalproblem1}--\eqref{originalproblem6}. Similarly, as the boundary $\partial\Omega_{u,v}$ is a curvilinear polygon of class $C^{1,1}$ (cf.\ \cite[Def.\ 1.4.5.1]{G85}) with the four vertices $$V=\{(-1,-1),(1,-1),(-1,0),(1,0)\}$$ connected by $W_\infty^2$-curves, it remains to check that $\omega_v\in(0,\pi/2]$, for all $v\in V$, where $\omega_v$ denotes the measure of the angle at $v\in V$, cf.\ \cite{WL12}. This is clearly satisfied as
\begin{align}
\omega_{(\pm 1,0)}  &=\arccos\left( \frac{(\partial_xu)(1)}{\sqrt{1+(\partial_xu)(1)^2}} \right)\in(0,\pi/2],\nonumber\\
\omega_{(\pm 1,-1)}&=\arccos\left(- \frac{(\partial_xv)(1)}{\sqrt{1+(\partial_xv)(1)^2}} \right)\in(0,\pi/2].\nonumber
\end{align}
It follows immediately that $u,v\in C^{2+\alpha}([-1,1])$ and, as $\Omega_{u,v}$ satisfies an exterior sphere condition (cf.\ \cite[p.\ 27]{GT}) on any boundary point, that $\phi\in C^{2+\alpha}(\Omega_{u,v}\cup\Gamma_u\cup\Gamma_v)$; cf.\ \cite[Thm.\ 6.13 \& Lem.\ 6.18]{GT}.
\endproof
\section{The small aspect ratio limit}
Fix $r_0$ and $\alpha$, $\lambda$ and $\mu$ as in Theorem~\ref{thm1} and consider a family $\{(u_\eps,v_\eps,\phi_\eps)\}_{\eps\in(0,1)}$ of solutions to \eqref{originalproblem1}--\eqref{originalproblem6}. Recall that $\tilde\phi_{\eps}=\phi_\eps\circ T^{-1}_{u_\eps,v_\eps}$, with $T_{u,v}$ as in \eqref{defT}, and define $\psi_\eps\colon\overline\Omega\to\R$ by $\psi_\eps(x',z')=\tilde\phi_{\eps}(x',z')-z'$. Let us first obtain some uniform estimates on $\psi_\eps$.
\lem\label{lem_estimates} Fix $0<\eps_0<1$. There exists a positive constant $K$ depending only on $\kappa_0$ and $\eps_0$ such that, for all $\eps\in(0,\eps_0]$,
\begin{align}
\norm{\psi_\eps}_{L_\infty(\Omega)}           & \leq 2,           \label{3.1}\\
\norm{\psi_\eps}_{L_2(\Omega)}                & \leq K\eps,       \label{3.2}\\
\norm{\partial_{z'}\psi_\eps}_{L_2(\Omega)}   & \leq K\eps,       \label{3.3}\\
\norm{\partial_{z'}^2\psi_\eps}_{L_2(\Omega)} & \leq K\eps^2.     \label{3.4}
\end{align}
\endlem
\proof By \eqref{estimate_tildephi}, we have $-2\leq\psi_\eps\leq 1$ and hence \eqref{3.1} follows. Let $f_\eps=f_{u_\eps,v_\eps;\eps}$ as in \eqref{deff}. Then
\beq\label{estimatef}\norm{f_\eps}_{L_\infty(\Omega)} \leq C_1\eps^2,\eeq
with $C_1=C_1(\kappa_0)>0$. In view of \eqref{Laplaceeq2}, \eqref{Laplacebc2} and \eqref{Lapdivform}, integration by parts yields
$$
\int_\Omega f_\eps\psi_\eps\d(x',z') = \int_\Omega(-\widetilde{\Delta}_\eps\psi_\eps)\psi_\eps\d(x',z') = \sum_{k=1}^5\mathcal I_k
$$
where
\begin{align}
\mathcal I_1 & = \eps^2\int_\Omega|\partial_{x'}\psi_\eps|^2\d(x',z'), \nonumber\\
\mathcal I_2 & = -2\eps^2\int_\Omega\frac{z'(\tilde u_\eps'-\tilde v_\eps')+\tilde v_\eps'}{\tilde u_\eps-\tilde v_\eps}(\partial_{x'}\psi_\eps)(\partial_{z'}\psi_\eps)\d(x',z'), \nonumber\\
\mathcal I_3 & = \int_\Omega\frac{1+\eps^2[z'(\tilde u_\eps'-\tilde v_\eps')+\tilde v_\eps']^2}{(\tilde u_\eps-\tilde v_\eps)^2}|\partial_{z'}\psi_\eps|^2\d(x',z'), \nonumber\\
\mathcal I_4 & = \eps^2 \int_\Omega(\tilde u_\eps'-\tilde v_\eps')\frac{z'(\tilde u_\eps'-\tilde v_\eps')+\tilde v_\eps'}{(\tilde u_\eps-\tilde v_\eps)^2}(\partial_{z'}\psi_\eps)\psi_\eps\d(x',z'), \nonumber\\
\mathcal I_5 & = -\eps^2 \int_\Omega\frac{\tilde u_\eps'-\tilde v_\eps'}{\tilde u_\eps-\tilde v_\eps}(\partial_{x'}\psi_\eps)\psi_\eps\d(x',z'). \nonumber
\end{align}
Now
$$
\mathcal I_1 + \mathcal I_2 + \mathcal I_3  = \int_\Omega\left[\eps^2\left(\partial_{x'}\psi_\eps-\frac{z'(\tilde u_\eps'-\tilde v_\eps')+\tilde v_\eps'}{\tilde u_\eps-\tilde v_\eps}\partial_{z'}\psi_\eps\right)^2
+\frac{|\partial_{z'}\psi_\eps|^2}{(\tilde u_\eps-\tilde v_\eps)^2}\right]d(x',z')
$$
and, as $|\psi_\eps|\leq 2$ and $|\Omega|= 2$,
$$|\mathcal I_4| \leq C_2\eps^2\norm{\partial_{z'}\psi_\eps}_{L_2(\Omega)},\quad C_2=C_2(\kappa_0)>0.$$
%
%
Again using integration by parts and \eqref{Laplacebc2}, $\mathcal I_5$ can be rewritten in the form
$$\mathcal I_5 = \frac{\eps^2}{2}\int_\Omega\partial_{x'}\left(\frac{\tilde u_\eps'-\tilde v_\eps'}{\tilde u_\eps-\tilde v_\eps}\right)\psi_\eps^2\d(x',z')$$
and hence there is a constant $C_3=C_3(\kappa_0)>0$ such that $|\mathcal I_5|\leq C_3\eps^2$. Now
\begin{align}
\int_\Omega f_\eps\psi_\eps\d(x',z')
& \geq \norm{\partial_{z'}\psi_\eps}_{L_2(\Omega)}^2
-C_2\eps^2\norm{\partial_{z'}\psi_\eps}_{L_2(\Omega)}-C_3\eps^2\nonumber\\
& = \left(\eps\norm{\partial_{z'}\psi_\eps}_{L_2(\Omega)}
-C_2\eps/2\right)^2+(1-\eps^2)\norm{\partial_{z'}\psi_\eps}_{L_2(\Omega)}^2\nonumber\\
& \qquad - \left(C_2^2/4+C_3\right)\eps^2\nonumber
\end{align}
and hence
$$4C_1\eps^2 \geq \int_\Omega f_\eps\psi_\eps\d(x',z') \geq
(1-\eps^2)\norm{\partial_{z'}\psi_\eps}_{L_2(\Omega)}^2 - \left(C_2^2/4+C_3\right)\eps^2.$$
As $0<\eps\leq\eps_0<1$ we infer that
$$\norm{\partial_{z'}\psi_\eps}_{L_2(\Omega)}^2\leq \frac{4C_1+C_2^2/4+C_3}{(1-\eps_0)^2}\eps^2.$$
Thus there is $C_4=C_4(\kappa_0,\eps_0)>0$ such that $\norm{\partial_{z'}\psi_\eps}_{L_2(\Omega)}\leq C_4\eps$ and this proves \eqref{3.3}. Estimate \eqref{3.2} follows from \eqref{Laplacebc2}, \eqref{3.3}, the Cauchy-Schwarz inequality and Fubini's Theorem as
\begin{align}
\norm{\psi_\eps}_{L_2(\Omega)}^2 & = \int_\Omega\left|\int_{z'}^1\partial_{z'}\psi_\eps(x',y')\d y'\right|^2d(x',z')\leq \int_{\Omega\times[0,1]}|\partial_{z'}\psi_\eps(x',y')|^2\d(x',y',z') \nonumber\\
& = \norm{\partial_{z'}\psi_\eps}_{L_2(\Omega)}^2\leq C_4^2\eps^2.\nonumber
\end{align}
Now we compute
\beq\label{eqmathcalJ}\int_\Omega f_\eps(\partial_{z'}^2\psi_\eps)\d(x',z') = \int_\Omega(-\widetilde{\Delta}_\eps\psi_\eps)(\partial_{z'}^2\psi_\eps)\d(x',z')=\sum_{k=1}^4\mathcal J_k\eeq
where
\begin{align}
\mathcal J_1 & = \eps^2\int_\Omega(\partial_{x'}^2\psi_\eps)(\partial_{z'}^2\psi_\eps)\d(x',z'), \nonumber\\
\mathcal J_2 & = -2\eps^2\int_\Omega\frac{z'(\tilde u_\eps'-\tilde v_\eps')+\tilde v_\eps'}{\tilde u_\eps-\tilde v_\eps}(\partial_{x'}\partial_{z'}\psi_\eps)(\partial_{z'}^2\psi_\eps)\d(x',z'), \nonumber\\
\mathcal J_3 & = \int_\Omega\frac{1+\eps^2(z'(\tilde u_\eps'-\tilde v_\eps')+\tilde v_\eps')^2}{(\tilde u_\eps-\tilde v_\eps)^2}(\partial_{z'}^2\psi_\eps)^2\d(x',z'), \nonumber\\
\mathcal J_4 & = \int_\Omega f_\eps(\partial_{z'}\psi_\eps)(\partial_{z'}^2\psi_\eps)\d(x',z'). \nonumber
\end{align}
By \cite[Lem.\ 4.3.1.2. \& 4.3.1.3]{G85}, $\mathcal J_1$ can be rewritten as
$$\mathcal J_1=\eps^2\int_\Omega (\partial_{x'}\partial_{z'}\psi_\eps)^2\d(x',z').$$
Hence \autoref{eqmathcalJ} can be rewritten as
\begin{align}
&\int_\Omega f_\eps(1-\partial_{z'}\psi_\eps)(\partial_{z'}^2\psi_\eps)\d(x',z') \nonumber\\
& \qquad = \eps^2\int_\Omega\left(\partial_{x'}\partial_{z'}\psi_\eps-\frac{z'(\tilde u_\eps'-\tilde v_\eps')+\tilde v_\eps'}{\tilde u_\eps-\tilde v_\eps}\partial_{z'}^2\psi_\eps\right)^2d(x',z')
+\int_\Omega\frac{|\partial_{z'}^2\psi_\eps|^2}{(\tilde u_\eps-\tilde v_\eps)^2}\d(x',z')\nonumber
\end{align}
and we conclude
\begin{align}
\norm{\partial_{z'}^2\psi_\eps}_{L_2(\Omega)}^2 & \leq \int_\Omega \frac{|\partial_{z'}^2\psi_\eps|^2}{(\tilde u_\eps-\tilde v_\eps)^2}\d(x',z')
\leq \int_\Omega f_\eps(1-\partial_{z'}\psi_\eps)(\partial_{z'}^2\psi_\eps)\d(x',z') \nonumber\\
& \leq \norm{f_\eps}_{L_\infty(\Omega)}(\sqrt 2+\norm{\partial_{z'}\psi_\eps}_{L_2(\Omega)})\norm{\partial_{z'}^2\psi_\eps}_{L_2(\Omega)}. \nonumber
\end{align}
By \eqref{3.3} and \eqref{estimatef}, there is a constant $C_5=C_5(\kappa_0,\eps_0)>0$ such that
$$\norm{\partial_{z'}^2\psi_\eps}_{L_2(\Omega)} \leq C_5\eps^2.$$
This implies \eqref{3.4} and choosing $K$ to be the maximum of $C_4$ and $C_5$ completes the proof of our lemma.
\endproof
\rem For the proof of Theorem~\ref{thm2} we will also need a slightly extended version of \cite[Lem.\ 12]{WL12}: For $\vartheta\in W_2^2(\Omega)$ we have
\beq
\norm{(\partial_{z'}\vartheta)(\cdot,0)}_{L_2(-1,1)},\norm{(\partial_{z'}\vartheta)(\cdot,1)}_{L_2(-1,1)}\leq\sqrt 2\left(\norm{\partial_{z'}\vartheta}_{L_2(\Omega)}+\norm{\partial_{z'}^2\vartheta}_{L_2(\Omega)}\right).\label{theta}
\eeq
Since the proof of the estimate for $\norm{(\partial_{z'}\vartheta)(\cdot,0)}_{L_2(-1,1)}$ is analogous to what is presented in \cite[Lem.\ 12]{WL12}, we omit it for the convenience of the reader.
\endrem
\emph{Proof of Theorem~\ref{thm2}:} Since $\{u_\eps\}_{\eps\in(0,1)}$ and $\{v_\eps\}_{\eps\in(0,1)}$ are bounded in $W_\infty^2(-1,1)$, a successive application of the Arzel\`{a}-Ascoli Theorem yields the existence of a null sequence $\{\eps_k\}_{k\geq 1}$, sequences $\{u_{\eps_k}\}_{k\geq 1}$ and $\{v_{\eps_k}\}_{k\geq 1}$ of real-valued functions on $(-1,1)$, and $u_0,v_0\in W_\infty^2(-1,1)$ such that
\begin{align}
u_{\eps_k} \to u_0 & \quad\text{in }W_\infty^1(-1,1), \nonumber\\
v_{\eps_k} \to v_0 & \quad\text{in }W_\infty^1(-1,1), \nonumber\\
u_{\eps_k} \stackrel{*}{\rightharpoonup} u_0 & \quad\text{in }W_\infty^2(-1,1). \nonumber\\
v_{\eps_k} \stackrel{*}{\rightharpoonup} v_0 & \quad\text{in }W_\infty^2(-1,1). \nonumber
\end{align}
Then $0\geq u_0(x)\geq -1/3+\kappa_0$ and $-1/3-2\kappa_0\geq v_0(x)\geq -1$ on $[-1,1]$ and, by Lemma~\ref{lem_estimates} and \eqref{theta},
$$
\frac{1+\eps_k^2(\partial_{x'}\tilde u_{\eps_k})^2}{(\tilde u_{\eps_k}-\tilde v_{\eps_k})^2}|(\partial_{z'}\tilde\phi_{\eps_k})(\cdot,1)|^2,
\frac{1+\eps_k^2(\partial_{x'}\tilde v_{\eps_k})^2}{(\tilde u_{\eps_k}-\tilde v_{\eps_k})^2}|(\partial_{z'}\tilde\phi_{\eps_k})(\cdot,0)|^2\to\frac{1}{(\tilde u_0-\tilde v_0)^2}
$$
in $L_1(-1,1)$. In particular, $u_0,v_0$ are solutions to the small aspect ratio equations \eqref{usar} and \eqref{vsar} satisfying $u_0(\pm 1)=0$ and $v_0(\pm 1)=-1$. Applying Lemma~\ref{lem_estimates} to $\psi_{\eps_k}=\tilde\phi_{\eps_k}-z'$, we get
$$\lim_{k\to\infty}\norm{\psi_{\eps_k}}_{L_2(\Omega)}^2=\lim_{k\to\infty}\int_\Omega|\tilde\phi_{\eps_k}(x',z')-z'|^2\d(x',z')=0$$
and by the change of variables $(x',z')\to(x,z)$
\begin{align}
\int_\Omega|\tilde\phi_{\eps_k}(x',z')-z'|^2\d(x',z') & = \int_{-1}^1\int_{v_{\eps_k}}^{u_{\eps_k}}\left|\phi_{\eps_k}(x,z)-\frac{z-v_{\eps_k}}{u_{\eps_k}-v_{\eps_k}}\right|^2\frac{d x\d z}{u_{\eps_k}-v_{\eps_k}}
\nonumber\\
& \geq \int_{-1}^1\int_{v_{\eps_k}}^{u_{\eps_k}}\left|\phi_{\eps_k}(x,z)-\frac{z-v_{\eps_k}}{u_{\eps_k}-v_{\eps_k}}\right|^2d x\d z,\nonumber
\end{align}
so that $\phi_{\eps_k}\mathbf{1}_{\Omega_{u_{\eps_k},v_{\eps_k}}}\to\phi_0\mathbf{1}_{\Omega_{u_0,v_0}}$ in $L_2(\Omega)$. This achieves the proof.\hfill$\square$
\section{Outlook}
In this paper we have shown that the problem \eqref{originalproblem1}--\eqref{originalproblem6} possesses solutions for small values $\lambda,\mu$ and we have given a justification of the associated small aspect ratio model by proving convergence towards solutions of \eqref{usar}, \eqref{vsar} when the aspect ratio of the device tends to zero. In this section, we present some plans for further research associated with the model introduced in this paper.

We do not yet have a deeper insight into the behavior of the MEMS model with two independent parameters $\lambda,\mu$ for large values of $\lambda$ or $\mu$. For instance, one can ask for the occurrence of the pull-in instability in case of, e.g., $\lambda$ small and $\mu$ large. A first goal for further work is to vary $(\lambda,\mu)$ through the parameter space $(0,\infty)^2$ and to check whether (and where) the solution constructed here ceases to exist and which consequences this has for the problem.

For the idealized MEMS model with a fixed ground plate, the authors of \cite{ELW12,ELW13} have recently shown that the parabolic free boundary problem with the condition
$$\partial_tu-\partial_x^2u = -\lambda\left(\eps^2|\partial_x\phi|^2+|\partial_z\phi|^2\right),\quad t>0,\quad x\in(-1,1),$$
on the free boundary possesses some interesting mathematical and physical properties: First, the associated MEMS model is shown to be locally well-posed in $\mathring{W}_p^2(-1,1)$, $2<p<\infty$. Next, the authors address global existence issues and present a result on asymptotic stability of the stationary solutions. Finally, the small aspect ratio limit is discussed in a similar manner as in \cite{WL12}. From the technical point of view, the approach of \cite{ELW12} benefits from a fixed point argument, regularizing effects of the heat semigroup, Lyapunov theory and the Principle of Linearized Stability. The parabolic MEMS model with one free boundary and an additional curvature term is the subject of \cite{ELW13}.

It is an open problem to study the parabolic free boundary problem with two free boundaries. Then the functions $u$ and $v$ depend on a temporal variable $t\geq 0$ and the spatial variable $x\in[-1,1]$ and the elliptic operator $-\partial_x^2$ on the left-hand side of \eqref{originalproblem3}, \eqref{originalproblem4} and \eqref{usar}, \eqref{vsar} has to be replaced by the parabolic operator $\partial_t-\partial_x^2$. A further generalization may be obtained by adding terms of the form $\alpha^2\partial_t^2u$ and $\beta\partial_x^4u$ on the left-hand side of \eqref{originalproblem3} and similarly in \eqref{originalproblem4} respectively, as explained in \cite{LL12}. The parameters $\alpha,\beta>0$ model the mechanical effects damping and bending.

In \cite{P01b,P01c}, a disk shaped idealized MEMS device model is presented. In this model, the rigid ground plate and the elastic membrane have drum shaped geometry. Assuming cylindrical symmetry, the displacement of the free membrane only depends on the radial component. It is not far to seek that the drum shaped model can be generalized by allowing for two moving boundaries as explained in this paper. A further reference concerning radially symmetric solutions of a MEMS model with bending and stretching is \cite{WL13} where the authors discuss fourth-order elliptic equations.\\[.25cm]
\textbf{Acknowledgement.} The author is grateful to the anonymous referee whose suggestions helped to improve the initial version of this paper.

\end{document}